\pdfoutput=1
\documentclass{mathprint}
\usepackage{macros}

\title[Uniformity of microsets]{On the uniformity and size of microsets}

\author
  {Richárd Balka}
  {HUN-REN Alfréd Rényi Institute of Mathematics, Reáltanoda u.~13--15, H-1053 Budapest, Hungary AND Institute of Mathematics and Informatics, Eszterházy Károly Catholic University, Leányka u.~4, H-3300 Eger, Hungary}
  {balkaricsi@gmail.com}

\author
  {Vilma Orgoványi}
  {
      Department of Stochastics,
      Institute of Mathematics,
      Budapest University of Technology and Economics,
      Műegyetem rkp.\ 3.,
      H-1111 Budapest,
      Hungary
  }
  {orgovanyi.vilma@gmail.com}

\author
  {Alex Rutar}
  {
      Department of Mathematics and Statistics,
      University of Jyväskylä,
      P.O.\ Box 35 (MaD),
      FI-40014 University of Jyväskylä,
      Finland
  }
  {alex@rutar.org}

\begin{document}
\begin{abstract}
    We resolve a few questions regarding the uniformity and size of microsets of subsets of Euclidean space.
    First, we construct a compact set $K\subset\R^d$ with Assouad dimension arbitrarily close to $d$ such that every microset of $K$ has no Ahlfors--David regular subset with dimension strictly larger than $0$.
    This answers a question of Orponen (also stated explicitly in \cite[Question~17.3.1]{zbl:1467.28001}).
    Then, we show that for any non-empty compact set $K\subset\R^d$ with lower dimension $\beta$, there is a microset $E$ of $K$ with finite $\beta$-dimensional packing pre-measure.
    This answers a strong version of a question of Fraser--Howroyd--Käenmäki--Yu \cite[Question~7.3]{zbl:1428.28013}, who previously obtained a similar result concerning the upper box dimension.
\end{abstract}

\section{Introduction}
One of the most fundamental geometric objects in analysis is the notion of a \emph{tangent}.
Generalizations, such as the \emph{weak tangent} and the \emph{microset}, have a long history in the study of the geometry of rough and irregular sets and metric spaces, with origins in the work of Furstenberg \cite{zbl:0208.32203} and Gromov \cite{zbl:0474.20018}.

More precisely, Furstenberg defined the $\ast$-dimension as the maximal Hausdorff dimension of a \emph{microset}, where a microset is any limit in the Hausdorff metric of \emph{minisets} of a compact set $K$, which are non-empty compact sets $\lambda(K-x)\cap B(0,1)$ for some $x\in K$ and $\lambda\geq 1$.
Here, $B(x,r)$ denotes the closed ball centred at $x$ with radius $r$.
We denote the set of microsets of $K$ by $\mathcal{G}_K$.
We refer to a microset for which the corresponding expansion ratios $\lambda$ diverge to infinity as a \emph{weak tangent}, and denote the set of weak tangents by $\Tan(K)$.
The related machinery of ``dynamical magnification'' has played a key role in the resolution of various long-standing conjectures in fractal geometry; for a (certainly incomplete) list of notable examples, see \cite{zbl:1251.28008,zbl:1409.11054,zbl:1426.11079,zbl:1430.11106}.

Parallel to the notion of a microset is the \emph{Assouad dimension}.
It is the value
\begin{equation*}
    \begin{aligned}
        \dimA K=\inf\Bigl\{\alpha>0&:(\exists C>0)\,(\forall 0<r\leq R<1)\,(\forall x\in K)\\
                                   &N_r\bigl(B(x,R)\cap K)\bigr)\leq C\Bigl(\frac{R}{r}\Bigr)^\alpha\Bigr\}
        \end{aligned}
\end{equation*}
where $N_r(F)$ is the number of closed balls of radius $r$ centred in $F$ required to cover $F$.
Heuristically, the Assouad dimension captures the worst-case scaling of $K$ at all locations and all scales.
In this form, the Assouad dimension has origins in embedding theory \cite{zbl:0396.46035} and was originally used to give a partial answer to the question concerning bi-Lipschitz embeddability of a metric space in Euclidean space.
Especially in recent years, the Assouad dimension has also become an important tool in quasiconformal analysis and fractal geometry.
Further references, along with many connections and applications to other fields, can be found in a variety of recent books on the topic \cite{zbl:1467.28001,zbl:1222.37004,zbl:1201.30002}.

The two concepts of Assouad dimension and microsets are intimately related.
The proof of the following result is due to Furstenberg \cite{zbl:1154.37322}, but the connection to Assouad dimension was made explicitly by Käenmäki, Ojala, \& Rossi.
\begin{proposition}[\cite{zbl:1154.37322,doi:10.1093/imrn/rnw336}]\label{p:Assouad-micro}
    Let $d\in\N$ and $K\subset\R^d$ be non-empty and compact.
    Then there is a weak tangent $E\in\Tan(K)$ such that $\mathcal{H}^{\dimA K}(E)>0$.
\end{proposition}
Here, $\mathcal{H}^s$ denotes the $s$-dimensional Hausdorff measure.
Since $\dimA E\leq\dimA K$ for any microset $E$, such a maximal microset necessarily has equal Assouad and Hausdorff dimensions.

This classification of Assouad dimension is very useful in practice, and plays an important role when studying the Assouad dimension.
To give one such example, Orponen \cite{zbl:1465.28008} proves the following result.
For $e\in S^1$, let $\pi_e\colon\R^2\to\R$ denote the orthogonal projection $\pi_e(x)=x\cdot e$.
Let $K\subset \mathbb{R}^2$, then the Hausdorff dimension of the set of directions for which the projection has smaller Assouad dimension than expected is 0:
\begin{equation*}
    \dimH\{e\in S^1:\dimA \pi_e(K)<\min\{\dimA K,1\}\}=0.
\end{equation*}
A key idea here is that the exceptional set bound is much stronger than the corresponding results with Hausdorff dimension in place of Assouad dimension (the sharp bound for Hausdorff dimension was recently established in \cite[Theorem~1.2]{arxiv:2308.08819}), primarily since in a heuristic sense one may reduce such a projection result to a projection result concerning microsets, which by \cref{p:Assouad-micro} enjoy much greater regularity than the original set.

In particular, it is of interest to understand precisely how much regularity we can expect microsets to have in general.
With \cref{p:Assouad-micro} in mind, we wish to answer the following question: to what extent can we guarantee additional homogeneity of a microset beyond that guaranteed by the Hausdorff dimension?

\subsection{Lower dimension and uniformity of microsets}
In order to quantify regularity more precisely, we begin by discussing an alternative notion of dimension which in some sense is the dual of Assouad dimension: the \emph{lower dimension}.
This notion of dimension was first introduced by Larman \cite{zbl:0152.24502}.
It is defined by
\begin{equation*}
    \begin{aligned}
        \dimL K=\sup\Bigl\{\alpha>0&:(\exists C>0)\,(\forall 0<r\leq R<1)\,(\forall x\in K)\\
                                   &N_r\bigl(B(x,R)\cap K)\bigr)\geq C\Bigl(\frac{R}{r}\Bigr)^\alpha\Bigr\}.
        \end{aligned}
\end{equation*}
Unfortunately, the lower dimension is not monotone under inclusion: in fact, if $K_1$ and $K_2$ are disjoint compact sets, then $\dimL (K_1\cup K_2)=\min\{\dimL K_1,\dimL K_2\}$.
 A natural way to rectify non-monotonicity under inclusion is simply to guarantee monotonicity by considering the \emph{modified lower dimension}, defined by
\begin{equation*}
    \dimML K=\sup\{\dimL E:\varnothing\neq E\subset K\}.
\end{equation*}
This notion of dimension first appeared in \cite{zbl:1390.28019}.
For compact sets $K$, it always holds that $\dimML K\leq\dimH K$ and moreover it is easy for the equality to be strict in general.
(Compactness is essential since the lower dimension is unchanged under taking the closure.)
After \cite{doi:10.1093/imrn/rnw336} appeared, Orponen asked whether or not one could necessarily find a weak tangent with an Ahlfors--David regular subset of the expected dimension.
A slightly weaker version of this question is also explicitly asked in \cite[Question~17.3.1]{zbl:1467.28001}.
Our first result states that the answer, in general, is no.
\begin{itheorem}\label{it:non-ahlfors}
    For any $d\in\N$ and $0\leq\alpha<d$, there exists a non-empty compact set $K\subset\R^d$ such that $\dimH K=\dimA K=\alpha$, but $\dimML E=0$ for all $E\in\mathcal{G}_K$.
    In particular, for all $s>0$ and $E\in\mathcal{G}_K$, $E$ does not contain an Ahlfors--David $s$-regular subset.
\end{itheorem}
The proof can be found in \cref{s:ad}.
The assumption that $\alpha<d$ is necessary: any subset of $\R^d$ with Assouad dimension $d$ has a microset with positive Lebesgue measure by \cref{p:Assouad-micro}, which in turn has a microset (in fact, a \emph{tangent}) with non-empty interior by the Lebesgue density theorem.
See also \cite[Theorem~2.4]{zbl:1429.11022} for a direct proof which avoids \cref{p:Assouad-micro} and the Lebesgue density theorem.

Next, we consider the analogous dual question concerning the relationship between microsets and lower dimensions.
In \cite{zbl:1428.28013}, the following result is shown.
\begin{proposition}[\cite{zbl:1428.28013}]
    Let $K\subset\R^d$ be a non-empty compact set.
    Then there exists an $F\in\Tan(K)$ such that $\dimuB F = \dimL K$.
\end{proposition}
Here, $\dimuB$ denotes the upper box (or Minkowski) dimension.
In \cite[Question~7.3]{zbl:1428.28013}, it is also asked whether there exists a weak tangent with finite $\dim_L(K)$-dimensional Hausdorff measure.
In fact, one might hope for more since the natural dual to positivity of the Hausdorff measure is finiteness of \emph{packing} measure.
This is our second result.
\begin{itheorem}\label{it:lower-tan}
    Let $K\subset\R^d$ be a non-empty compact set with lower dimension $\beta$.
    Then there exists an $F\in\Tan(K)$ such that $\mathcal{P}^\beta(F)<\infty$.
\end{itheorem}
Actually, we obtain an upper bound for the packing pre-measure depending only on $\beta$; see \cref{t:pack}.
In particular, the packing pre-measure and measure are the same by \cite{zbl:1018.28004}.
Proofs of these results, as well as precise definitions of packing (pre-)measures, can be found in \cref{ss:lower-micro}.

In the usual proof of the corresponding result for Assouad dimension (that is, \cref{p:Assouad-micro}), one constructs a certain measure satisfying a Frostman-type property on the dyadic tree representing $K$, and then takes a pushforward of that measure onto the set itself.
The main technical difficulty is that such an approach does not immediately work for packing dimension, since it could happen that the set $K$ is badly aligned with respect to the dyadic tree, and such a push-forward could be much smaller than expected.
Instead, we use a generalized system of cubes, similar to those first introduced by David \cite{zbl:0696.42011} and Christ \cite{zbl:0758.42009}.
We use this tree representation to identify a set of locations and scales with desirable properties, and then construct an appropriate tangent measure directly on the set $K$.

\section{Microsets and Ahlfors--David regular subsets}\label{s:ad}
In this section, we prove our first result, \cref{it:non-ahlfors}.

\subsection{Dyadic microsets}\label{ss:dyadic}
We begin by introducing a standard representation of a compact set using dyadic cubes.
Fix $d\in\N$ and a non-empty compact set $K\subset[0,1]^d\subset\R^d$.
Let $\mathcal{D}=\bigcup_{n=0}^\infty\mathcal{D}_n$ denote the set of closed dyadic cubes which intersect $K$, where $\mathcal{D}_n$ denotes the subset of dyadic cubes with side-length $2^{-n}$ intersecting $K$.
Since each dyadic cube in $\mathcal{D}_n$ can contain at most $2^d$ dyadic cubes in $\mathcal{D}_{n+1}$, there is a natural coding of $\mathcal{D}$ using words in $\{1,\ldots,2^d\}^*$ which gives $\mathcal{D}$ the structure of a tree.

Given a dyadic cube $Q\in\mathcal{D}$, let $T_Q$ denote the homothety satisfying $T_Q(Q) = [0,1]^d$.
\begin{definition}
    We say that a compact set $E\subset[0,1]^d$ is a \emph{dyadic microset} if it is a limit in the Hausdorff metric of sets $T_{Q_{n}}(K)\cap [0,1]^d$ for $Q_{n}\in\mathcal{D}$.
\end{definition}
We make two standard observations which essentially mean, in this section, that we need only consider dyadic microsets and dyadic covering numbers.
\begin{enumerate}[nl]
    \item Since each ball $B(x,r)$ can be covered by a finite number of dyadic cubes of side length approximately $r$, by a pigeonholing argument it follows that the Assouad dimension can be rephrased dyadically: for a non-empty compact set $K\subset[0,1]^d$, $\dimA K$ is the infimum over all constants $\alpha\geq 0$ so that for all $k\in\N\cup\{0\}$, dyadic cubes $Q\in\mathcal{D}_k$, and all $m\in\N\cup\{0\}$,
        \begin{equation*}
            \#\{R\in \mathcal{D}_{k+m}:R\subset Q\}\leq C_{\alpha} 2^{m\alpha},
        \end{equation*}
        where $C_{\alpha}$ is a constant depending only on $\alpha$.
    \item For any microset $E\in\mathcal{G}_K$, there is a homothety $h(x)=\lambda x+ t$ with $1\leq \lambda<2$ and a set of dyadic microsets $E_1,\ldots,E_k$ where $k\leq 4^d$ with translations $t_1,\ldots,t_k\in\R^d$ so that
        \begin{equation*}
            h(E)\subset (E_1+t_1)\cup\cdots\cup(E_k+t_k).
        \end{equation*}
\end{enumerate}

\subsection{Sets with small microsets}\label{ss:non-ahlfors}
We are now ready to prove \cref{it:non-ahlfors}.

We will construct a ``uniformly branching set'' which has scaling properties which are inhomogeneous in scale, but homogeneous in space.
Such constructions are widespread in the literature, under names such as ``homogeneous Moran sets'' or ``non-autonomous self-similar sets''.
We begin by presenting this construction in general.
Let $\bm{a}\in\{0,1\}^{\N}$: we define a compact set $K(\bm{a})$ as follows.
We will define a nested sequence of compact sets $(K_n)_{n=0}^\infty$ where each $K_n$ is a union of dyadic sub-cubes of $[0,1]^d$ at level $n$.
Begin with $K_0=[0,1]^d$.
Then for each $n\in\N\cup\{0\}$ and level-$n$ dyadic cube $Q$ of $K_n$, if $a_{n+1}=0$, replace $Q$ with a single dyadic cube $Q'\subset Q$ at level $n+1$, and if $a_{n+1}=1$, replace $Q$ with $2^d$ dyadic cubes at level $n+1$.
(The specific choice of the sub-cube $Q'\subset Q$ is not important.)
Finally, let $K=\bigcap_{n=0}^\infty K_n$.
Equivalently, the dyadic tree corresponding to $K$ has full branching at each level $n$ where $a_n=1$, and no branching where $a_n=0$.

We separate the construction into two parts: in \cref{sss:branching} we construct the branching numbers $\bm{a}\in\{0,1\}^{\N}$ and in \cref{sss:properties} we show that the corresponding set has the desired properties.

\subsubsection{Constructing the branching numbers}\label{sss:branching}
Given a sequence $\bm{x}=(x_n)_{n=1}^\infty$, we denote the \emph{lower Cesàro mean} by
\begin{equation*}
    \underline{\lambda}(\bm{x}) =\liminf_{n\to\infty}\frac{\sum_{j=1}^n x_{j}}{n}.
\end{equation*}
We begin by constructing a sequence $\bm{a}=(a_n)_{n=1}^\infty\in\{0,1\}^{\N}$ with the following properties.
\begin{lemma}\label{l:seq-construct}
    For every $\varepsilon>0$, there exists a binary sequence $\bm{a}=(a_n)_{n=1}^\infty\in\{0,1\}^{\N}$ and a sequence $(N_m)_{m=1}^\infty$ of natural numbers such that:
    \begin{enumerate}[nl,r]
        \item\label{im:eps} $\underline{\lambda}(\bm{a})\geq 1-\varepsilon$, and
        \item\label{im:consec} For all $m\in\N$ and $j\in\N$, the sequence $(a_{j},\ldots,a_{j+N_m-1})$ contains a consecutive sequence of $m$ zeros.
    \end{enumerate}
\end{lemma}
\begin{proof}
    First, let $\gamma\in(0,1)$ and for each $m\in\N$, let
    \begin{equation*}
        N_m=\lfloor m^3\gamma^{-1}\rfloor.
    \end{equation*}
    Note that $N_m\geq m$ since $\gamma<1$.
    Now define sequences $\bm{a}_m=(a_{n,m})_{n=1}^\infty$ where $\bm{a}_m=(\nu_m,\nu_m,\nu_m,\ldots)$ and each $\nu_m$ is a finite sequence of $N_m-m$ 1s followed by $m$ 0s.
    Observe moreover that
    \begin{equation*}
        \underline{\lambda}(\bm{a}_m) \geq 1-\frac{m}{N_m}\geq 1-\frac{\gamma}{m^2}.
    \end{equation*}
    Now define the sequence $\bm{a}=(a_n)_{n=1}^\infty$ by
    \begin{equation*}
        a_n =
        \begin{cases}
            0 &: a_{n,m}=0\text{ for some }m\in\N,\\
            1 &:\text{ otherwise.}
        \end{cases}
    \end{equation*}
    Clearly the sequence $\bm{a}$ satisfies \cref{im:consec}.
    Moreover,
    \begin{equation*}
        \underline{\lambda}(\bm{a})\geq 1-\sum_{m=1}^\infty\bigl(1-\underline{\lambda}(\bm{a}_m)\bigr)\geq 1-\sum_{m=1}^\infty\frac{\gamma}{m^2}\geq 1-\frac{\gamma\pi^2}{6},
    \end{equation*}
    which, for $\gamma$ sufficiently small, can be made arbitrarily close to $1$.
\end{proof}

\subsubsection{Completing the example}\label{sss:properties}
It remains to prove that the set $K(\bm{a})$ associated with a sequence satisfying the conclusions of \cref{l:seq-construct} has the desired properties.
\begin{theorem}\label{t:small-microset}
    For every $\alpha\in[0,d)$, there exists a set $K\subset\R^d$ with $\dimA K=\alpha$ such that for every microset $E\in\mathcal{G}_K$, $\dimML E=0$.
\end{theorem}
\begin{proof}
    If $\alpha=0$ this result is trivial, so we may fix $0<\alpha<d$.
    By \cref{l:seq-construct}, get a sequence $\bm{a}$ and natural numbers $(N_m)_{m=1}^\infty$ such that $\underline{\lambda}(\bm{a})\geq \alpha/d$ and for each $m\in\N$ and $j\in\N$, the sequence $(a_j,\ldots,a_{j+N_m-1})$ contains a consecutive sequence of $m$ zeros.
    Let $K=K(\bm{a})$ denote the corresponding uniformly branching set.
    At level $n$, $K$ intersects $2^{d \sum_{i=1}^n a_i}$ distinct dyadic cubes, so
    \begin{equation*}
        \dimA K\geq d\cdot \liminf_{n\to\infty}\frac{\sum_{i=1}^n a_i}{n}=d\cdot\underline{\lambda}(\bm{a})\geq \alpha.
    \end{equation*}

    Now fix an arbitrary microset $E\in\mathcal{G}_K$.
    Since the modified lower dimension is finitely stable and monotonic under inclusion, we may assume that $E$ is a dyadic microset.
    Fix a subset $\varnothing\neq F\subset E$: it suffices to prove that $\dimL F=0$.

    Write $E=\lim_{n\to\infty}K_{Q_n}$ for $Q_n\in\mathcal{D}$.
    Let $m\in\N$ be arbitrary, and get $k=k(m)$ so that
    \begin{equation}\label{e:dh-close}
        d_{\mathcal{H}}(K_{Q_k}, E)\leq 2^{-N_m}.
    \end{equation}
    Let $j$ be such that $Q_k\in\mathcal{D}_j$.
    By construction of $\bm{a}$, there is an $i\in\Z$ with
    \begin{equation*}
        j\leq i+1\leq i+m\leq j+N_m-1
    \end{equation*}
    so that $a_{i+1}=\cdots=a_{i+m}=0$.

    Now fix an arbitrary $x_0\in F$ and consider the subset $B\bigl(x_0,2^{j-i}\bigr)\cap F$.
    Let $y_0=T_{Q_k}(x_0)$, so that $B(y_0,2^{-i})$ can be covered by $3^d$ dyadic cubes in level $i$.
    Fix one such dyadic cube $Q$.
    Since $a_{i+1}=\cdots=a_{i+m}=0$, $Q\cap K$ can be covered by 1 dyadic cube at level $i+m$.
    Since $i+m\leq j+N_m-1$, the $2^{-(j+N_m)}$-neighbourhood of $Q\cap K$ can thus be covered by $3^d$ dyadic cubes at level $i+m$.
    Thus by \cref{e:dh-close}, since $F\subset E$ and each dyadic cube at level $i+m$ is contained in a ball of radius $2^{-(i+m)}$,
    \begin{equation*}
        N_{2^{j-i+m}}\left(B\bigl(x_0,2^{j-i}\bigr)\cap F\right)=N_{2^{-i+m}}\left(B(y_0,2^{-i})\cap T_{Q_k}^{-1}(E)\right)\leq 9^d.
    \end{equation*}
    But $m\in\N$ was arbitrary, so that $\dimL F=0$, as required.

    In the above, we constructed a set $K(\bm{a})$ satisfying the desired properties such that $\dimA K(\bm{a})\geq\alpha$.
    To construct a set $K$ with $\dimA K=\alpha$, we perform the same construction as above, but instead by subdividing each cube into either $1$ or $2^d$ sub-cubes (depending on the value of $a_n$) with side-lengths $\rho\in(0,1/2]$.
    Call the resulting set $K_\rho(\bm{a})$.
    Then, standard arguments (see, for example, \cite[Lemma~3.2]{zbl:1371.28024} or \cite[Theorem~2]{zbl:1405.28007}) show that
    \begin{equation*}
        \dimA K_\rho(\bm{a})=\frac{d\log(2)}{\log(1/\rho)}\cdot\limsup_{n\to\infty}\sup_{k\in\N\cup\{0\}}\frac{a_{k+1}+\cdots+a_{k+n}}{n}.
    \end{equation*}
    In particular, $\dimA K_\rho(\bm{a})$ is a continuous and decreasing function of $\rho$, so there exists a value $\rho_0\in(0,1/2]$ such that the corresponding set $K=K_{\rho_0}(\bm{a})$ has $\dimA K=\alpha$.
    Moreover, the above proof works identically with the modified construction (using $\rho_0$ in place of $1/2$ where appropriate), giving the desired result.
\end{proof}
\begin{remark}
    In the above proof, we have essentially shown that $\dimML K(\bm{a}) = \dimL K(\bm{a})$ and moreover, that every microset of $K(\bm{a})$ is essentially of the form $K(\bm{b})$ where $\bm{b}$ is a limit of finite sub-sequences of $\bm{a}$.
\end{remark}
We can now conclude the proof of \cref{it:non-ahlfors}.
\begin{proofref}{it:non-ahlfors}
    Fix $\alpha\in(0,d)$.
    First by \cref{t:small-microset}, for a set $K\subset\R^d$ with $\dimA K=\alpha$ such that for each $E\in\mathcal{G}_K$, $\dimML E=0$.
    By \cref{p:Assouad-micro}, there is a microset $F\in\mathcal{G}_K$ such that $\dimA F=\dimH F=\alpha$.
    But an easy computation shows that every microset of $F$ is a subset of a microset of $K$ (for an explicit proof, see \cite[Lemma~3.11]{zbl:1342.28016}).
    In particular, $\dimML E=0$ for all $E\in\mathcal{G}_F$.
    Finally, since any Ahlfors--David $s$-regular set has lower dimension $s$, every $E\in\mathcal{G}_F$ and $s$-regular subset of $E$ must have $s=0$.
\end{proofref}

\section{Small microsets for packing measure}\label{ss:lower-micro}
In this section, we prove \cref{it:lower-tan}.

\subsection{Preliminaries on packing dimension}

Let $s\geq 0$ be fixed.
We first define \emph{$s$-dimensional packing pre-measure} for an arbitrary set $E\subset\R^d$ as
\begin{equation*}
    \mathcal{P}_0^s(E)=\lim_{\delta\to 0}\sup\left\{\sum_{i=1}^\infty (2r_i)^s:\begin{matrix}\{B(x_i,r_i)\}\text{ pairwise disjoint}\\
    \text{with $r_i\leq\delta$ and $x_i\in E$}\end{matrix}\right\}.
\end{equation*}
Note that the limit always exists by monotonicity.
In general, if $\mathcal{P}^s_0(E)<\infty$, then $\dimuB E\leq s$ and if $\mathcal{P}^s_0(E)>0$, then $\dimuB E \geq s$.
Countably stabilizing the pre-measure yields the \emph{$s$-dimensional packing measure}:
\begin{equation*}
    \mathcal{P}^s(K)=\inf\left\{\sum_{i=1}^\infty\mathcal{P}_0^s(E_i):K\subset\bigcup_{i=1}^\infty E_i\right\}.
\end{equation*}
For completeness, we also recall the definition of the \emph{packing dimension}:
\begin{equation*}
    \dimP K=\inf\{s\geq 0:\mathcal{P}^s(K)=0\}.
\end{equation*}
Finally, we recall the following simple `dual' Frostman property for packing pre-measure.
\begin{lemma}\label{l:anti-frostman}
    Let $K\subset\R^d$ be arbitrary and let $c>0$ and $r_0\in(0,1)$ be fixed.
    Suppose $\mu$ is a Borel probability measure satisfying $\mu(B(x,r)) \geq c r^s$ for all $x\in K$ and $r\in(0,r_0)$.
    Then
    \begin{equation*}
        \mathcal{P}^s_0(K) \leq c^{-1} 2^s
    \end{equation*}
\end{lemma}
\begin{proof}
    Let $\{B(x_i,r_i)\}_{i=1}^\infty$ be an arbitrary packing of $K$ where $0<r_i<r_0$.
    Then
    \begin{align*}
        1 = \mu(\R^d)\geq \mu\left(\bigcup_{i=1}^\infty B(x_i,r_i)\right) = \sum_{i=1}^\infty \mu(B(x_i,r_i))\geq c 2^{-s} \sum_{i=1}^\infty (2r_i)^s.
    \end{align*}
    Since the packing was arbitrary, it follows that $\mathcal{P}^s_0(K) \leq 2^s c^{-1}$.
\end{proof}

\subsection{Notation for trees}
We now introduce some basic notation concerning trees.

Fix some integer $d\in\N$ and a number $\rho\in(0,1)$.
Write $\mathcal{A} = \{1,\ldots,d\}$, and set $\mathcal{A}^*=\bigcup_{n=0}^\infty\mathcal{A}^n$.
We equip the space $\Delta=\{1,\ldots,d\}^{\N}$ with the metric
\begin{equation*}
    d(x,y) = \inf\{ \rho^n: x_j=y_j\text{ for all }j=1,\ldots,n\}.
\end{equation*}
This metric makes $\Delta$ into a compact ultrametric space, with topology corresponding to the usual product topology on $\Delta$.
For $\mtt{a} = (a_1,\ldots,a_n)\in\mathcal{A}^n$, we denote the corresponding open and closed cylinder
\begin{equation*}
    [\mtt{a}] \coloneqq B(x,\rho^n) = \left\{x \in\Delta: x_j = a_j\text{ for }j=1,\ldots,n\right\}.
\end{equation*}
Given a subset $\mathcal{T}\subset\mathcal{A}^*$, we write $\mathcal{T}_n = \mathcal{T}\cap\mathcal{A}^n$.

Now suppose $\Omega$ is a non-empty compact subset of $\Delta$.
We associate with $\Omega$ a tree $\mathcal{T}=\bigcup_{n=0}^\infty\mathcal{T}_n$, which is just the collection of cylinders
\begin{equation*}
    \mathcal{T}=\left\{[\mtt{a}]:\mtt{a}\in \mathcal{A}^*;\, [\mtt{a}]\cap \Omega\neq\varnothing\right\}.
\end{equation*}
Conversely, given a non-empty subset $\mathcal{T}\subset\mathcal{A}^*$ with the property that each $\mtt{a}\in\mathcal{T}_n$ has at least one suffix $\mtt{a}j\in\mathcal{T}_{n+1}$, the \emph{boundary} of $\mathcal{T}$ is the set
\begin{equation*}
    \partial\mathcal{T} = \bigcap_{n=0}^\infty\bigcup_{\mtt{a}\in\mathcal{T}_n} [\mtt{a}].
\end{equation*}
Clearly the tree associated with $\partial\mathcal{T}$ is precisely $\mathcal{T}$.

Given $Q=[\mtt{a}]$ where $\mtt{a}\in\mathcal{T}$, let:
\begin{itemize}[nl]
    \item $M_Q\colon Q\to\Delta$ denote the magnification defined by $M_Q(\mtt{a} x) = x$.
    \item $\mathcal{T}(Q)\coloneqq\{\mtt{b}\in\mathcal{T}:[\mtt{b}]\subset Q\}$, and
    \item $\mathcal{T}^Q\coloneqq\{M_Q(Q'):Q'\in\mathcal{T}(Q)\}$.
\end{itemize}
We call a set $\mathcal{T}^Q$ for $Q\in\mathcal{T}$ a \emph{subtree}.

\subsection{Generalized cubes for compact sets}
In \cref{ss:dyadic}, we defined the dyadic tree associated with a compact set $K$.
However, dyadic cubes can be aligned badly with respect to a set $K$, so microsets on the corresponding dyadic tree $\mathcal{D}$ do not correspond in a reasonable way with microsets of $K$.

In order to work around this problem, we use a system of generalized cubes, an idea first introduced by David \cite{zbl:0696.42011} and Christ \cite{zbl:0758.42009}.
Let us begin with a definition highlighting the properties of such a representation that we require.
Unlike the usual properties of such trees, we require no hypotheses concerning the measures of the boundaries of the sets.
The most important property here which is not satisfied by the usual system of dyadic cubes is the first inclusion in \cref{e:incl}.
\begin{definition}
    Let $K\subset\R^d$ be a non-empty compact set and let $0<\rho<1$.
    Let $\mathcal{B}=\{(Q_{i,k},x_{i,k}):k\in\N\cup\{0\};\,i\in \mathcal{N}_k;\, x_{i,k}\in Q_{i,k}\}$ be a family of non-empty Borel sets with distinguished points $x_{i,k}\in K$.
    We say that $\mathcal{B}$ is an \emph{inner regular partition} if there are constants $c,C>0$ so that the following hold:
    \begin{enumerate}[nl,r]
        \item $\#\mathcal{N}_0=1$.
        \item For all $k\in\N\cup\{0\}$, $K=\bigcup_{i\in\mathcal{N}_k}Q_{i,k}$ disjointly.
        \item The system of Borel sets has a \emph{tree structure}: if $k\leq m$, $i\in\mathcal{N}_k$, and $j\in\mathcal{N}_m$, then either $Q_{i,k}\cap Q_{j,m}=\varnothing$ or $Q_{i,k}\supset Q_{j,m}$.
        \item For all $k\in\N\cup\{0\}$ and $i\in\mathcal{N}_k$,
            \begin{equation}\label{e:incl}
                B(x_{i,k}, c\rho^k)\subset Q_{i,k} \subset B(x_{i,k}, C\rho^k).
            \end{equation}
        \item We have $\{x_{i,k}:i\in\mathcal{N}_k\}\subset\{x_{i,k+1}:i\in\mathcal{N}_{k+1}\}$.
    \end{enumerate}
\end{definition}
That such a tree structure exists for arbitrary compact sets can be found explicitly in \cite[Theorem~2.1]{zbl:1277.28017}.
Tracking constants and applying an initial rescaling yields the following.
\begin{proposition}\label{p:tree-repr}
    Let $K\subset\R^d$ be a non-empty compact set and let $\rho=1/4$.
    Suppose $\diam K \leq 1$.
    Then $K$ has an inner regular partition $\mathcal{B}$ for any constants $0<c<1/6$ and $C\geq 4/3$.
\end{proposition}
Moreover, since our set $K$ is a subset of a doubling metric space, we note the following straightforward bound on the branching numbers.
\begin{lemma}\label{l:branch-bound}
    Suppose $\mathcal{B}$ is an inner regular partition for $K$.
    Then there is a constant $M\in\N$ so that for all $k\in\N\cup\{0\}$ and $i\in\mathcal{N}_k$,
    \begin{equation*}
        \#\{j\in\mathcal{N}_{k+1}:Q_{j,k+1}\subset Q_{i,k}\} \leq M.
    \end{equation*}
\end{lemma}
Now, let $\mathcal{B}$ be an inner regular partition of $K$ with respect to a constant $\rho\in(0,1)$, and let $M$ be as in \cref{l:branch-bound}.
The structure $\mathcal{B}$ induces the structure of a tree on $K$ as follows.

We begin by fixing a labelling of the cylinders $Q_{i,k}$.
Set $\mathcal{A} = \{1,\ldots,M\}$.
First, assign to the unique root cylinder $Q_0$ the word $\varnothing$.
Now, inductively, recalling \cref{l:branch-bound}, having assigned words $\mtt{a}(Q)\in\mathcal{A}^k$ for all $Q=Q_{i,k}$, if $Q_{i,k}$ has children $R_{1},\ldots,R_j$ with $j\leq M$, set $\mtt{a}(R_j) = \mtt{a}(Q) j$.
Clearly this assignment is a bijection, so given $k\in\N$ and a word $\mtt{a}\in\mathcal{A}^k$, we may speak of the corresponding cylinder $Q(\mtt{a})$ and distinguished point $x(\mtt{a})\in Q(\mtt{a})$.
Let $\mathcal{T}$ denote the corresponding tree, with boundary $\partial\mathcal{T}$.
We call such tree a $(\rho, M)$-tree, which we will use without making the set $K$ explicit.
We define the lower dimension of a boundary of a tree analogously to sets:
\begin{equation*}
    \begin{aligned}
        \dimL \partial \mathcal{T} =\sup\Bigl\{s>0&:(\exists C>0)\,(\forall 0\leq m\leq k)\,(\forall [\mtt{a}]\in \mathcal{T}_m)\\
                                   &\#\bigl( \mtt{b}\in \mathcal{T}_k, [\mtt{b}]\subset [\mtt{a}] \bigr)\geq C\Bigl(\rho^{(m-k)}\Bigr)^s\Bigr\}.
        \end{aligned}
\end{equation*}

To conclude this section, we note the following lemma.
\begin{lemma}\label{l:lower-tree-equiv}
    Let $K\subset\R^d$ be a non-empty compact set with inner regular partition $\mathcal{B}$ and corresponding tree $\mathcal{T}$.
    Then $\dimL\partial\mathcal{T} \leq \dimL K$.
\end{lemma}
\begin{proof}
    Let $\mathcal{B}$ be an inner regular partition with respect to the parameters $\rho\in(0,1)$ and constants $0<c\leq C$.
    Let $s<\dimL\partial\mathcal{T}$ be arbitrary.
    By definition of the lower dimension, there is a constant $\delta>0$ so that for all integers $0\leq m\leq k$ and $\mtt{a}\in\mathcal{T}_m$,
    \begin{equation}\label{e:tree-lower}
        \#\{\mtt{b}\in \mathcal{T}_k: [\mtt{b}]\subset [\mtt{a}]\} \geq \delta \rho^{(m-k)s}.
    \end{equation}
    Set
    \begin{equation*}
        0<\eta =\min\left\{\left(\frac{c\rho}{2C}\right)^s, \delta \left(\frac{ c\rho^2}{2C}\right)^s\right\},
    \end{equation*}
    and let $y\in K$ and $0<r\leq R<1$ be arbitrary.
    We will show that
    \begin{equation*}
        N_r(B(y,R)\cap K) \geq \eta\left(\frac{R}{r}\right)^s.
    \end{equation*}
    First note that if $r\geq c(2C)^{-1}\rho R$, then
    \begin{equation*}
        \eta\left(\frac{R}{r}\right)^s \leq 1 \leq N_r(B(y,R)\cap K).
    \end{equation*}
    Thus, it remains to handle the case that $r< c(2C)^{-1}\rho R$.
    Let $m\in\N$ be minimal so that $2C\rho^m \leq R$, and let $\mtt{a}\in\mathcal{T}_m$ be chosen so that $y\in Q(\mtt{a})$.
    In particular, since $Q(\mtt{a})\subset B(x(\mtt{a}), C\rho^m)$, it follows that $Q(\mtt{a})\subset B(y,R)$.

    Now, the assumption on $r$ ensures $r< c\rho^m$.
    Let $k\geq m$ be maximal so that $r<c\rho^k$.
    Also recall that $[\mtt{a}] = B(z,\rho^m)$ for any $z\in\partial\mathcal{T}\cap [\mtt{a}]$.
    Thus, by \cref{e:tree-lower}, get $\{\mtt{b}_1,\ldots,\mtt{b}_M\}\subset\mathcal{T}_k$ such that $[\mtt{b}_j]\subset [\mtt{a}]$ and $M\geq \delta \rho^{(m-k)s}$.
    Moreover, for each $\mtt{b}_j$, observe that
    \begin{equation*}
        B(x(\mtt{b}_j), r)\cap K\subset B(x(\mtt{b}_j), c\rho^k)\cap K \subset Q(\mtt{b}_j) \subset Q(\mtt{a})\subset B(y,R)\cap K.
    \end{equation*}
    Therefore since the balls $B(x(\mtt{b}_j),c\rho^k)$ are disjoint,
    \begin{equation*}
        N_r(B(y,R)\cap K) \geq M \geq \delta \rho^{(m-k)s} \geq \delta\left(\frac{\rho R}{2C}\cdot \frac{c \rho}{r}\right)^s \geq \eta\left(\frac{R}{r}\right)^s,
    \end{equation*}
    as required.
\end{proof}

\subsection{Pigeonholing good scales and locations}
In order to obtain the upper bound on the packing pre-measure, we would like to use \cref{l:anti-frostman}.
However, it will not actually be useful for us to construct a measure on $\partial\mathcal{T}$ itself.
We will simply use the tree $\mathcal{T}$ to inform a good choice of cylinders.
This lemma is the dual of Furstenberg's construction of Frostman measures on microsets; see for instance the exposition in \cite[Lemma 2.4.4]{zbl:1390.28012}.
\begin{lemma}\label{l:reverse-furstenberg}
    Let $\mathcal{T}$ be a $(\rho, M)$-tree, and let $\ell\in\N$ be fixed.
    Let $0<t<s$ be arbitrary, and suppose $k\geq \ell t/(s-t)$ is such that
    \begin{equation}\label{e:top-count}
        \#\mathcal{T}_k \leq \rho^{-k s}.
    \end{equation}
    Then, there are $0\leq n \leq k-\ell$ and $Q\in\mathcal{T}_n$ such that for all $j=0,\ldots,\ell$ and $Q'\in(\mathcal{T}^Q)_j$,
    \begin{equation}\label{e:antif}
        \frac{\#(Q'\cap \mathcal{T}^Q_{k-n})}{\#\mathcal{T}^Q_{k-n}} \geq \rho^{t j}.
    \end{equation}
\end{lemma}
\begin{proof}
    Set $Q_0=[0,1]^d$.
    If $Q_0$ satisfies \cref{e:antif}, we are done.
    Otherwise, there is an $1\leq\ell_1\leq \ell$ and a $Q_1\in\mathcal{T}_{\ell_1}$ so that \cref{e:antif} fails.
    Equivalently, since $\#\mathcal{T}^{Q_0}_{k} = \# T^{Q_1}_{k-\ell_1}$,
    \begin{equation}\label{e:sub-count}
        \#\mathcal{T}^{Q_1}_{k-\ell_1}< \rho^{-k s}\rho^{\ell_1 t}.
    \end{equation}
    Repeating the above procedure with the iterated bound \cref{e:sub-count} in place of \cref{e:top-count}, either there is some $n\leq k-\ell$ and a cylinder $Q_m\in\mathcal{T}_n$ which satisfies \cref{e:antif}, or $k\geq\ell_1+\cdots+\ell_m>k-\ell$.
    Suppose for contradiction that the latter situation occurs.
    Then
    \begin{equation*}
        1 \leq \#\mathcal{T}^{Q_m}_{k-(\ell_1+\cdots+\ell_m)}<\rho^{-k s}\rho^{(\ell_1+\cdots+\ell_m) t}\leq \rho^{-ks}\rho^{(k-\ell) t}.
    \end{equation*}
    Rearranging, $k < \ell t/(s-t)$ which contradicts the choice of $k$.
\end{proof}
Applying the definition of the lower dimension yields the following.
\begin{corollary}\label{c:good-measure}
    For all $t>\dimL\partial\mathcal{T}$ and $\ell\in\N$, there are $n\geq \ell$, $Q\in\mathcal{T}_n$, and $k\geq n+\ell$ such that for all $j=0,\ldots,\ell$ and $Q'\in\mathcal{T}^Q_j$,
    \begin{equation*}
        \frac{\#(Q'\cap \mathcal{T}^Q_{k-n})}{\#\mathcal{T}^Q_{k-n}} \geq \rho^{t j}.
    \end{equation*}
\end{corollary}
\begin{proof}
    Let $t>\dimL\partial\mathcal{T}$ and $\ell\in\N$, and fix $\dimL\partial\mathcal{T}< s < t$.
    Let $k_0 \geq \ell t/(s-t)$ be fixed.
    Since $s>\dimL\partial\mathcal{T}$, for all $\delta>0$, there exists integers $0\leq n_0\leq m$ and $P_0\in\mathcal{T}_{n_0}$ so that
    \begin{equation*}
        1 \leq \#\mathcal{T}^{P_0}_{m-n_0}<\delta \rho^{(n_0-m)s}.
    \end{equation*}
    In particular, taking $\delta = \rho^{(k_0+\ell)s}$, we must have $m-n_0\geq k_0+\ell$.
    Now let $n_1 = n_0+\ell$ and let $P\in\mathcal{T}^{P_0}_{\ell}$ be arbitrary.
    Then
    \begin{equation*}
        \#\mathcal{T}^P_{m - n_1} < \delta \rho^{(n_0-m)s} \leq \delta \rho^{-\ell s} \rho^{(n_1-m)s} < \rho^{(n_1-m)s}.
    \end{equation*}
    Crucially, observe that $m-n_1\geq k_0$.
    Therefore, we may apply \cref{l:reverse-furstenberg} to the tree $\mathcal{T}^P$ to obtain a cylinder $Q_0\in\mathcal{T}^P$ satisfying \cref{e:antif}.
    But if $P$ has coding $\mtt{a}$ and $Q_0$ has coding $\mtt{b}$, let $Q=[\mtt{a}\mtt{b}]$ so $\mathcal{T}^Q =(\mathcal{T}^P)^{Q_0}$.
    Then since $P\in\mathcal{T}_n$ where $n\geq \ell$, $Q$ has the desired properties.
\end{proof}

\subsection{Constructing good measures}
Finally, we can complete the proof of \cref{it:lower-tan}.
\begin{theorem}\label{t:pack}
    Let $K\subset\R^d$ be non-empty and compact, with $\beta=\dimL K$.
    Then there exists an $F\in\Tan(K)$ such that
    \begin{equation*}
        \mathcal{P}^\beta_0(F)\leq 257^\beta<\infty.
    \end{equation*}
\end{theorem}
\begin{proof}
    First, observe that $\Tan(K)$ is unchanged upon applying a homothety to $K$.
    Therefore, we may assume that $\diam K \leq 1$.
    Fix $\rho=1/4$.
    Applying \cref{p:tree-repr}, get an inner regular partition $\mathcal{B}$ with corresponding constants $c$ and $C$.
    Let $\mathcal{T}$ be the corresponding tree, so that $\dimL\partial\mathcal{T}\leq\beta$ by \cref{l:lower-tree-equiv}.
    In particular, by \cref{c:good-measure}, for all $\ell\in\N$, there is an $n(\ell)\geq \ell$, $Q_\ell\in\mathcal{T}_{n(\ell)}$, and $m(\ell)\coloneqq k-n(\ell) \geq\ell$ such that for all $j=0,\ldots,\ell$ and $Q'\in\mathcal{T}^{Q_\ell}_j$,
    \begin{equation*}
        \frac{\#(Q'\cap \mathcal{T}^{Q_\ell}_{m(\ell)})}{\#\mathcal{T}^{Q_{\ell}}_{m(\ell)}} \geq \rho^{(\beta+1/\ell)j}.
    \end{equation*}
    Let $x_\ell\in Q_\ell$ be such that
    \begin{equation*}
        B(x_\ell, c \rho^{n(\ell)})\subset Q_\ell\subset B(x_\ell, C\rho^{n(\ell)}).
    \end{equation*}
    Next, denote the set of centres at level $m(\ell)$ contained in $Q_\ell$ by
    \begin{equation*}
        \mathfrak{C}_\ell = \{x_{i,m(\ell)}: i\in\mathcal{N}_{m(\ell)};\,x_{i,m(\ell)}\in Q_\ell\}.
    \end{equation*}
    We then define the measure
    \begin{equation*}
        \mu_\ell = \frac{1}{\#\mathfrak{C}_\ell} \sum_{z\in\mathfrak{C}_\ell}\delta_{z}.
    \end{equation*}
    Observe that $\mu_\ell$ is a probability measure, and $\supp \mu_\ell \subset Q_\ell$.
    By construction,
    \begin{equation*}
        \mu_\ell(Q') = \frac{\#(Q'\cap \mathcal{T}^Q_{m(\ell)})}{\#\mathcal{T}^Q_{m(\ell)}} \geq \rho^{(\beta+1/\ell)j}
    \end{equation*}
    for all $j=0,\ldots,\ell$ and $Q'\in(\mathcal{T}^{Q_\ell})_j$.
    Let $f_\ell$ be the unique homothety mapping $B(x_\ell, c \rho^{n(\ell)}/2)$ to $B(0,1)$, and set
    \begin{equation*}
        \nu_\ell = \mu_\ell \circ f_{\ell}^{-1}.
    \end{equation*}
    Since $\supp\mu_\ell\subset Q_\ell$, it follows that $\supp \nu_\ell \subset B(0, 2C c^{-1})$.
    Passing to a subsequence if necessary, we may set
    \begin{equation*}
        F = \lim_{\ell \to\infty} f_\ell(K)\cap B(0,1).
    \end{equation*}
    Here, the limit is with respect to the Hausdorff metric.
    Since $x_\ell\in K$ and $n(\ell)$ diverges to infinity, $F\in\Tan(K)$.
    Passing again to a subsequence, we may assume that $\lim_{\ell\to\infty}\nu_\ell = \nu$ in the weak-$*$ topology.

    Now, let $x\in F$ and $0<r<1$ be arbitrary.
    Let $N$ be sufficiently large so that for all $\ell\geq N$,
    \begin{enumerate}[nl,a]
        \item there is a $y_\ell\in f_\ell(K)\cap B(0,1)$ such that $d(y_\ell, x) \leq r/2$, and
        \item\label{im:rel} $\rho^N \leq cr / (4C)$.
    \end{enumerate}
    Let $z_\ell = f_\ell^{-1}(y_\ell) \in Q_\ell$ and let $j\in\N\cup\{0\}$ be minimal such that
    \begin{equation*}
        C \rho^{n(\ell)+j} \leq \frac{r\cdot c \rho^{n(\ell)}}{4}.
    \end{equation*}
    By the choice of $N$ in \cref{im:rel} and the minimality of $j$, it follows that $j\leq \ell$.
    Then let $Q\in\mathcal{T}_{n(\ell)+j}$ be the unique cylinder which contains $z_\ell$; of course, $Q\subset Q_\ell$.
    By construction of $\mu_\ell$ and since $j\leq \ell$, we have that
    \begin{equation*}
        \mu_\ell(Q) \geq \rho^{(\beta+1/\ell) j} \geq r^{\beta + 1/\ell} \cdot \left(\frac{c \rho}{4C}\right)^{\beta + 1/\ell}.
    \end{equation*}
    Moreover, since $f_\ell$ is a homothety with expansion factor $2 (c \rho^{n(\ell)})^{-1}$, $f_\ell(Q)$ has diameter at most $r/2$, so in fact $f_\ell(Q)\subset B(x,r)$.
    Therefore,
    \begin{equation*}
        \nu_{\ell}(B(x,r))\geq \nu_\ell(f_{\ell}(Q)) \geq r^{\beta + 1/\ell} \cdot \left(\frac{c \rho}{4C}\right)^{\beta + 1/\ell}.
    \end{equation*}
    It follows that
    \begin{equation*}
        \nu(B(x, r)) \geq \limsup_{\ell\to\infty}\nu_\ell(B(x,r)) \geq \left(\frac{c \rho}{4C}\right)^\beta r^\beta .
    \end{equation*}
    Since $x\in F$ and $0<r<1$ were arbitrary and $c$, $C$, and $\rho$ are fixed, we conclude by \cref{l:anti-frostman} that
    \begin{equation*}
        \mathcal{P}^\beta(F) \leq \mathcal{P}^\beta_0(F) \leq \left(\frac{8C}{c\rho}\right)^\beta<\infty.
    \end{equation*}
    Substituting the relevant constants from \cref{p:tree-repr} yields the result.
\end{proof}

\begin{acknowledgements}
    The authors thank Tuomas Orponen for interesting discussions regarding the topics in this paper.

    RB is supported by the National Research, Development and Innovation Office -- NKFIH, grant no.~146922 and the János Bolyai Research Scholarship of the Hungarian Academy of Sciences.
    VO is supported by the grant KKP144059 ``Fractal geometry and applications'', the grant NKFI K142169, and the Hungarian Research Network through the HUN-REN--BME Stochastics Research Group.
    AR is supported by the Research Council of Finland via Tuomas Orponen's project \emph{Approximate incidence geometry}, grant no.\ 355453.
\end{acknowledgements}

@preprint{arxiv:2308.08819,
  arxiv = {2308.08819},
  author = {Ren, Kevin and Wang, Hong},
  eprint = {2308.08819},
  eprinttype = {arxiv},
  month = {08},
  title = {Furstenberg sets estimate in the plane},
  year = {2023},
}

@article{doi:10.1093/imrn/rnw336,
  author = {Käenmäki, Antti and Ojala, Tuomo and Rossi, Eino},
  doi = {10.1093/imrn/rnw336},
  eprint = {10.1093/imrn/rnw336},
  eprinttype = {doi},
  journal = {Int. Math. Res. Not.},
  month = {6},
  pages = {3769--3799},
  title = {Rigidity of quasisymmetric mappings on self-affine carpets},
  volume = {2018},
  year = {2018},
}

@article{zbl:0758.42009,
  author = {Christ, Michael},
  doi = {10.4064/cm-60-61-2-601-628},
  eprint = {0758.42009},
  eprinttype = {zbl},
  journal = {Colloq. Math.},
  language = {English},
  pages = {601--628},
  title = {A {{\(T(b)\)}} theorem with remarks on analytic capacity and the {Cauchy} integral},
  volume = {60/61},
  year = {1990},
  zbl = {0758.42009},
  zbmath = {00010412},
}

@article{zbl:1018.28004,
  author = {Feng, De-Jun and Hua, Su and Wen, Zhi-Ying},
  doi = {10.1112/S0024609399006256},
  eprint = {1018.28004},
  eprinttype = {zbl},
  issue = {6},
  journal = {Bull. Lond. Math. Soc.},
  language = {English},
  pages = {665--670},
  publisher = {John Wiley \ Sons, Chichester; London Mathematical Society, London},
  title = {Some relations between packing premeasure and packing measure},
  volume = {31},
  year = {1999},
  zbl = {1018.28004},
  zbmath = {01463645},
}

@article{zbl:0152.24502,
  author = {Larman, D. G.},
  doi = {10.1112/plms/s3-17.1.178},
  eprint = {0152.24502},
  eprinttype = {zbl},
  journal = {Proc. Lond. Math. Soc.},
  language = {English},
  pages = {178--192},
  title = {A new theory of dimension},
  volume = {17},
  year = {1967},
  zbl = {0152.24502},
  zbmath = {03244600},
}

@inproceedings{zbl:0208.32203,
  author = {Furstenberg, Harry},
  booktitle = {Problems in Analysis},
  eprint = {0208.32203},
  eprinttype = {zbl},
  language = {English},
  number = {31},
  pages = {41--59},
  publisher = {Princeton University Press},
  series = {Princeton Math. Ser.},
  title = {Intersections of {Cantor} sets and transversality of semi-groups},
  year = {1970},
  zbl = {0208.32203},
  zbmath = {03331653},
}

@thesis{zbl:0396.46035,
  address = {Orsay},
  author = {Assouad, Patrice},
  eprint = {0396.46035},
  eprinttype = {zbl},
  institution = {Univ. Paris XI},
  language = {French},
  series = {Publ. Math. Orsay},
  title = {Espaces métriques, plongements, facteurs},
  type = {Thèse de doctorat d’État},
  year = {1977},
  zbl = {0396.46035},
  zbmath = {03615337},
}

@article{zbl:0474.20018,
  author = {Gromov, Mikhael},
  doi = {10.1007/BF02698687},
  eprint = {0474.20018},
  eprinttype = {zbl},
  journal = {Publ. Math., Inst. Hautes Étud. Sci.},
  language = {English},
  pages = {53--78},
  title = {Groups of polynomial growth and expanding maps},
  titleaddon = {Appendix by {Jacques} {Tits}},
  volume = {53},
  year = {1981},
  zbl = {0474.20018},
  zbmath = {03743523},
}

@article{zbl:0696.42011,
  author = {David, Guy},
  doi = {10.4171/RMI/64},
  eprint = {0696.42011},
  eprinttype = {zbl},
  journal = {Rev. Mat. Iberoam.},
  language = {French},
  pages = {73--114},
  title = {Morceaux de graphes lipschitziens et intégrales singulières sur une surface},
  volume = {4},
  year = {1988},
  zbl = {0696.42011},
  zbmath = {04140480},
}

@article{zbl:1154.37322,
  author = {Furstenberg, Hillel},
  doi = {10.1017/S0143385708000084},
  eprint = {1154.37322},
  eprinttype = {zbl},
  journal = {Ergodic Theory Dyn. Syst.},
  language = {English},
  pages = {405--422},
  title = {Ergodic fractal measures and dimension conservation},
  volume = {28},
  year = {2008},
  zbl = {1154.37322},
  zbmath = {05272907},
}

@book{zbl:1201.30002,
  address = {Providence, RI},
  author = {Mackay, John M. and Tyson, Jeremy T.},
  eprint = {1201.30002},
  eprinttype = {zbl},
  language = {English},
  number = {54},
  publisher = {American Mathematical Society},
  series = {Univ. Lect. Ser.},
  subtitle = {Theory and application},
  title = {Conformal dimension},
  year = {2010},
  zbl = {1201.30002},
  zbmath = {05731476},
}

@book{zbl:1222.37004,
  address = {Cambridge},
  author = {Robinson, James C.},
  eprint = {1222.37004},
  eprinttype = {zbl},
  language = {English},
  number = {186},
  publisher = {Cambridge University Press},
  series = {Camb. Tracts Math.},
  title = {Dimensions, embeddings, and attractors},
  year = {2011},
  zbl = {1222.37004},
  zbmath = {05835385},
}

@article{zbl:1251.28008,
  author = {Hochman, Michael and Shmerkin, Pablo},
  doi = {10.4007/annals.2012.175.3.1},
  eprint = {1251.28008},
  eprinttype = {zbl},
  journal = {Ann. Math.},
  language = {English},
  pages = {1001--1059},
  title = {Local entropy averages and projections of fractal measures},
  volume = {175},
  year = {2012},
  zbl = {1251.28008},
  zbmath = {06051266},
}

@article{zbl:1277.28017,
  author = {Käenmäki, Antti and Rajala, Tapio and Suomala, Ville},
  doi = {10.1090/S0002-9939-2012-11161-X},
  eprint = {1277.28017},
  eprinttype = {zbl},
  journal = {Proc. Am. Math. Soc.},
  language = {English},
  pages = {3275--3281},
  title = {Existence of doubling measures via generalised nested cubes},
  volume = {140},
  year = {2012},
  zbl = {1277.28017},
  zbmath = {06203928},
}

@article{zbl:1409.11054,
  author = {Hochman, Michael and Shmerkin, Pablo},
  doi = {10.1007/s00222-014-0573-5},
  eprint = {1409.11054},
  eprinttype = {zbl},
  journal = {Invent. Math.},
  language = {English},
  pages = {427--479},
  title = {Equidistribution from fractal measures},
  volume = {202},
  year = {2015},
  zbl = {1409.11054},
  zbmath = {06504096},
}

@article{zbl:1342.28016,
  author = {Käenmäki, Antti and Rossi, Eino},
  doi = {10.5186/aasfm.2016.4133},
  eprint = {1342.28016},
  eprinttype = {zbl},
  journal = {Ann. Acad. Sci. Fenn., Math.},
  language = {English},
  pages = {465--490},
  title = {Weak separation condition, {Assouad} dimension, and {Furstenberg} homogeneity},
  volume = {41},
  year = {2016},
  zbl = {1342.28016},
  zbmath = {06551802},
}

@book{zbl:1390.28012,
  address = {Cambridge},
  author = {Bishop, Christopher J. and Peres, Yuval},
  doi = {10.1017/9781316460238},
  eprint = {1390.28012},
  eprinttype = {zbl},
  language = {English},
  number = {162},
  publisher = {Cambridge University Press},
  series = {Camb. Stud. Adv. Math.},
  title = {Fractals in probability and analysis},
  year = {2017},
  zbl = {1390.28012},
  zbmath = {06653783},
}

@article{zbl:1371.28024,
  author = {Olson, Eric J. and Robinson, James C. and Sharples, Nicholas},
  doi = {10.1017/S0305004115000584},
  eprint = {1371.28024},
  eprinttype = {zbl},
  journal = {Math. Proc. Camb. Philos. Soc.},
  language = {English},
  pages = {51--75},
  title = {Generalised {Cantor} sets and the dimension of products},
  volume = {160},
  year = {2016},
  zbl = {1371.28024},
  zbmath = {06782156},
}

@article{zbl:1429.11022,
  author = {Fraser, Jonathan M. and Yu, Han},
  doi = {10.1112/blms.12112},
  eprint = {1429.11022},
  eprinttype = {zbl},
  journal = {Bull. Lond. Math. Soc.},
  language = {English},
  pages = {85--95},
  title = {Arithmetic patches, weak tangents, and dimension},
  volume = {50},
  year = {2018},
  zbl = {1429.11022},
  zbmath = {06846521},
}

@article{zbl:1390.28019,
  author = {Fraser, Jonathan M. and Yu, Han},
  doi = {10.1016/j.aim.2017.12.019},
  eprint = {1390.28019},
  eprinttype = {zbl},
  journal = {Adv. Math.},
  language = {English},
  pages = {273--328},
  title = {New dimension spectra: finer information on scaling and homogeneity},
  volume = {329},
  year = {2018},
  zbl = {1390.28019},
  zbmath = {06863445},
}

@article{zbl:1405.28007,
  author = {García, Ignacio and Hare, Kathryn and Mendivil, Franklin},
  doi = {10.1017/S0308210517000488},
  eprint = {1405.28007},
  eprinttype = {zbl},
  journal = {Proc. R. Soc. Edinb., Sect. A, Math.},
  language = {English},
  pages = {517--540},
  title = {Assouad dimensions of complementary sets},
  volume = {148},
  year = {2018},
  zbl = {1405.28007},
  zbmath = {07000762},
}

@article{zbl:1426.11079,
  author = {Shmerkin, Pablo},
  doi = {10.4007/annals.2019.189.2.1},
  eprint = {1426.11079},
  eprinttype = {zbl},
  journal = {Ann. Math.},
  language = {English},
  pages = {319--391},
  title = {On {Furstenberg}'s intersection conjecture, self-similar measures, and the {{\(L^q\)}} norms of convolutions},
  volume = {189},
  year = {2019},
  zbl = {1426.11079},
  zbmath = {07041748},
}

@article{zbl:1430.11106,
  author = {Wu, Meng},
  doi = {10.4007/annals.2019.189.3.2},
  eprint = {1430.11106},
  eprinttype = {zbl},
  journal = {Ann. Math.},
  language = {English},
  pages = {707--751},
  title = {A proof of {Furstenberg}'s conjecture on the intersections of {{\(\times p\)}}- and {{\(\times q\)}}-invariant sets},
  volume = {189},
  year = {2019},
  zbl = {1430.11106},
  zbmath = {07097489},
}

@article{zbl:1428.28013,
  author = {Fraser, Jonathan M. and Howroyd, Douglas C. and Käenmäki, Antti and Yu, Han},
  doi = {10.1090/proc/14613},
  eprint = {1428.28013},
  eprinttype = {zbl},
  journal = {Proc. Am. Math. Soc.},
  language = {English},
  pages = {4921--4936},
  title = {On the {Hausdorff} dimension of microsets},
  volume = {147},
  year = {2019},
  zbl = {1428.28013},
  zbmath = {07113905},
}

@book{zbl:1467.28001,
  address = {Cambridge},
  author = {Fraser, Jonathan M.},
  doi = {10.1017/9781108778459},
  eprint = {1467.28001},
  eprinttype = {zbl},
  language = {English},
  number = {222},
  publisher = {Cambridge University Press},
  series = {Camb. Tracts Math.},
  title = {Assouad dimension and fractal geometry},
  year = {2020},
  zbl = {1467.28001},
  zbmath = {07274738},
}

@article{zbl:1465.28008,
  author = {Orponen, Tuomas},
  doi = {10.1112/plms.12317},
  eprint = {1465.28008},
  eprinttype = {zbl},
  journal = {Proc. Lond. Math. Soc.},
  language = {English},
  pages = {317--351},
  title = {On the {Assouad} dimension of projections},
  volume = {122},
  year = {2021},
  zbl = {1465.28008},
  zbmath = {07348179},
}
\end{document}